\documentclass[a4paper,notitlepage,twoside,reqno,11pt]{amsart}

\usepackage{anysize}
\marginsize{3.4cm}{3.4cm}{3cm}{3cm}

\usepackage{bbm,pifont,latexsym,dcolumn,indentfirst}
\usepackage[hypertexnames=false,bookmarksopen=true,linktocpage=true,pdfstartview={XYZ null null 1.25}]{hyperref}
\usepackage{amsthm,amsmath,amssymb,amscd,amsfonts,mathrsfs}
\usepackage{color,graphicx,xcolor,graphics,subfigure,extarrows,caption2}
\usepackage{pdfpages,titletoc}

\usepackage{autonum}

\newtheorem{thm}{Theorem}[section]
\newtheorem{lem}[thm]{Lemma}
\newtheorem{cor}[thm]{Corollary}

\newtheorem*{mainthm}{Main Theorem}

\theoremstyle{definition}
\newtheorem*{defi}{Definition}

\newtheorem*{rmk}{Remark}

\newcommand{\Z}{\mathbb{Z}}
\newcommand{\Q}{\mathbb{Q}}
\newcommand{\R}{\mathbb{R}}
\newcommand{\D}{\mathbb{D}}

\newcommand{\C}{\mathbb{C}}
\newcommand{\EC}{\widehat{\mathbb{C}}}

\newcommand{\MC}{\mathcal{C}}
\newcommand{\MD}{\mathcal{D}}

\newcommand{\MI}{\mathcal{I}}
\newcommand{\MM}{\mathcal{M}}
\newcommand{\MO}{\mathcal{O}}
\newcommand{\MP}{\mathcal{P}}

\newcommand{\MS}{\mathcal{S}}
\newcommand{\MU}{\mathcal{U}}

\newcommand{\ii}{\textup{i}}
\newcommand{\id}{\textup{id}}

\newcommand{\Int}{\textup{int}}
\newcommand{\Ext}{\textup{ext}}

\newcommand{\Comp}{\textup{Comp}}

\newcommand{\bess}{\begin{eqnarray*}}
\newcommand{\eess}{\end{eqnarray*}}

\makeatletter\@addtoreset{equation}{section}\makeatother

\begin{document}

\author[X. Wang]{Xiaoguang Wang}
\address{School of Mathematical Sciences, Zhejiang University, Hangzhou 310027, P. R. China}
\email{wxg688@163.com}

\author[F. Yang]{FEI YANG}
\address{School of Mathematics, Nanjing University, Nanjing 210093, P. R. China}
\email{yangfei@nju.edu.cn}

\title[Capture components of cubic Siegel polynomials]{Capture components of cubic Siegel polynomials}

\begin{abstract}
Let $\theta$ be an irrational number of bounded type. We prove that all capture components in the parameter space of cubic polynomials $f_a(z)=e^{2\pi\ii\theta}z+a z^2+z^3$, where $a\in\C$, are Jordan domains.
\end{abstract}



\date{\today}


 \dedicatory{Dedicated to Professor Mitsuhiro Shishikura on the occasion of his $60$th birthday}

\maketitle


\section{Introduction}\label{introduction}

\subsection{Backgrounds}

The study of parameter spaces in the field of complex dynamics is an important topic. After Douady and Hubbard's pioneer work on quadratic polynomials $Q_c(z)=z^2+c$ with $c\in\C$ \cite{DH8485}, the Mandelbrot set
\begin{equation}
\mathbb{M}:=\{c\in\C: \text{the iteration } Q_c^{\circ n}(0)\not\to\infty \text{ as } n\to+\infty\}
\end{equation}
attracted many people.
For example, Shishikura proved that the boundary of the Mandelbrot set has Hausdorff dimension two \cite{Shi98}. McMullen proved that the Mandelbrot set is universal in the holomorphic families of rational maps \cite{McM00b}.
One of the most important conjectures in this field is that the Mandelbrot set is locally connected, which implies the density of hyperbolicity conjecture. See the work of Dudko and Lyubich \cite{DL23a}, \cite{DL26} and the references therein for the progress on these conjectures.

For high dimensional parameter spaces, Branner and Hubbard first gave a characterization of the global topology of the parameter space of cubic polynomials \cite{BH88}.
They decomposed the parameter space (complex $2$-dimension) into the connectedness locus and the union of leaves marked by the escaping rate of the critical points. For more topological characterizations of the connectedness locus of cubic polynomials, see \cite{Lav89}, \cite{Mil92}, \cite{EY99} and \cite{BOPT14}.

In general it is difficult to study the high dimensional complex parameter spaces since $\C^n$ is at least real $4$-dimensional if $n\geqslant 2$.
Some results in this direction include \cite{Mil12}, \cite{WY17}, \cite{CWY22}, where the topological structures of the high dimensional hyperbolic components are studied.
One way to understand high dimensional parameter spaces is to study their various one-dimensional slices.
In these slices, the bifurcation loci and the connected components of their complements are important study objects. In particular, the boundaries of the  complex one-dimensional hyperbolic components have been studied extensively. See \cite{Roe07}, \cite{QRWY15}, \cite{RWY17}, \cite{Wan21}, \cite{QRW23} and the references therein.

Especially for cubic polynomials, the related results are very affluent.
For example, one can find such meaningful slices in the parameter space of cubic polynomials with real coefficients \cite{Mil92}, with an irrationally indifferent fixed point \cite{Zak99}, \cite{BH01}, \cite{BOST22}, with a parabolic fixed point \cite{Nak05}, \cite{Lom14}, \cite{Zha24R}, and with a super-attracting cycle (see \cite{Fau92}, \cite{Roe07}, \cite{Mil09}, \cite{BKM10}, \cite{Wan21}). See also \cite{BOT22} for the slice of cubic polynomials with a non-repelling fixed point.

\subsection{Main results}

Let $f$ be a rational map with $\deg(f)\geqslant 2$ having an irrationally indifferent fixed point at the origin, i.e., $f(0)=0$ and $f'(0)=e^{2\pi\ii\theta}$ with $\theta\in\R\setminus\Q$.
If $f$ is conjugate to the rigid rotation $R_\theta(z)=e^{2\pi\ii\theta}z$ in a neighborhood of $0$, then $f$ is called locally linearizable at $0$, and the maximal region in which $f$ is conjugate to $R_\theta$ is a simply connected domain $\Delta(f)$ called the \textit{Siegel disk} of $f$ centered at $0$.

Let $\theta\in(0,1)\setminus\Q$ be an irrational number of \textit{bounded type}, i.e., the continued fraction expansion $[0;a_1,a_2,\cdots, a_n,\cdots]$ of $\theta$ satisfies $\sup_{n\geqslant 1}\{a_n\}<\infty$. According to \cite{Sie42}, every cubic polynomial in the slice
\begin{equation}\label{equ:f-a}
\MS(\theta):=\Big\{f_a(z)=e^{2\pi\ii\theta}z+a z^2+z^3:a\in\C\Big\}
\end{equation}
has a Siegel disk $\Delta(f_a)$ centered at $0$.
By \cite{Zak99}, the boundary of every $\Delta(f_a)$ is a quasi-circle passing through one or both critical points of $f_a$ (see \cite{Shi01} and \cite{Zha11} for more general results).
It may happen that the forward orbit of one critical point of $f_a$ enters into the Siegel disk $\Delta(f_a)$ eventually. Such parameters form a subset of $\MS(\theta)$ and they belong to the complement of the bifurcation locus although they are not hyperbolic.
Each component of these parameters is called a \textit{capture component} (see Figure \ref{Fig:parameter}).

\begin{figure}[!htpb]
  \setlength{\unitlength}{1mm}
  \centering
  \includegraphics[width=0.95\textwidth]{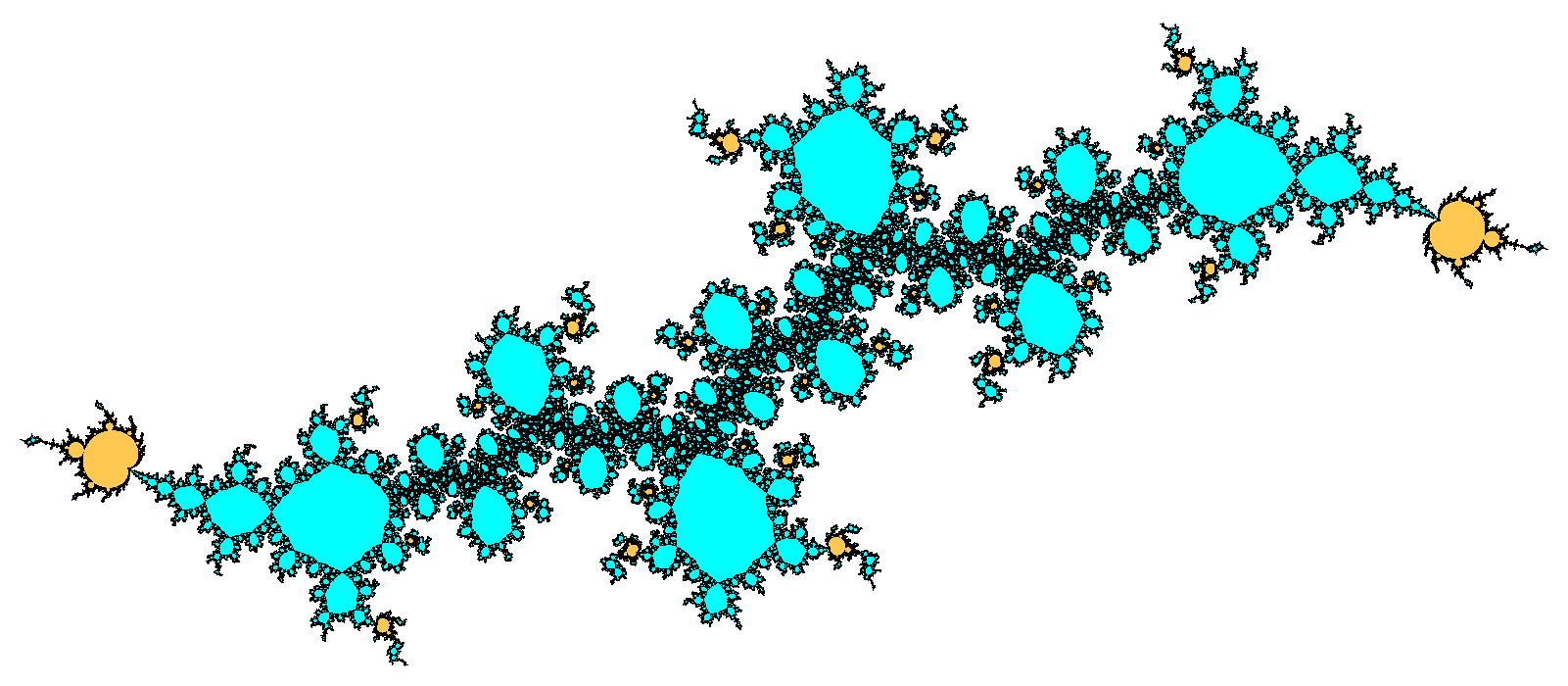}
  \caption{The parameter plane of $\MS(\theta)$ with $\theta=(\sqrt{5}-1)/2$. The capture components are colored cyan and the hyperbolic-like parameters (those $a$'s such that $f_a$ has a bounded attracting cycle) are colored yellow. This picture has been rotated.}
  \label{Fig:para-BH}
\end{figure}

In this paper, we study the topology of the boundaries of capture components of $\MS(\theta)$ for fixed bounded type $\theta$ and obtain the following result.

\begin{mainthm}\label{thm:main}
For any bounded type irrational number $\theta$, all capture components of $\MS(\theta)$ are Jordan domains.
\end{mainthm}

The parameter slice $\MS(\theta)$ has attracted many people.
Zakeri was the first who studied the connectedness locus of $\MS(\theta)$ when $\theta$ is of bounded type \cite{Zak99}.
Buff and Henriksen proved that for any real number $\theta$, the bifurcation locus of $\MS(\theta)$ contains a quasi-conformal copy of the Julia set of the quadratic polynomial $Q_\theta(z)=e^{2\pi\ii\theta}z+z^2$ \cite{BH01}.
This implies that some capture components of $\MS(\theta)$ are Jordan domains if the bounded Fatou components of $Q_\theta$ are.
This includes the case when $\theta$ is a rational number and an irrational number but of bounded type, Petersen-Zakeri type \cite{PZ04} and sufficiently high type \cite{Che25}, \cite{SY25}.

By constructing puzzles, J. Yang and R. Zhang proved that for bounded type $\theta$, the locus of non-renormalizable maps in $\MS(\theta)$ is homeomorphic to a double-copy of the filled Julia set of $Q_\theta$ (minus the Siegel disk) glued along the Siegel boundary \cite{YZ24}, which verifies a conjecture of Blokh-Oversteegen-Ptacek-Timorin \cite{BOPT14}.
In particular, this implies that some capture components of $\MS(\theta)$ are Jordan domains, and one may also obtain the Main Theorem by combining \cite{BH01} and \cite{YZ24}. However, the proof in this paper does not rely on puzzles.

Recently, for $\theta\in\Q$, R. Zhang proved that all capture components in $\MS(\theta)$ are Jordan domains \cite{Zha24R}.
For more results on the characterizations of the connectedness locus of $\MS(\theta)$, see \cite{Zak18}, \cite{Che20B}, \cite{BCOT21} and the references therein.

\medskip
Based on the work of Douady-Herman (see \cite{Dou87}, \cite{Her87}) and the Main Theorem, a natural question is: are all capture components of $\MS(\theta)$ quasi-disks for bounded type $\theta$? We prove that the answer is affirmative for those capture components whose closures are disjoint with the Zakeri curve (i.e., the parameters such that the corresponding Siegel disk boundaries contain two critical points). See \S\ref{sec:quasi-disk}. For the capture components whose closures intersect with the Zakeri curve, the holomorphic motion argument cannot be used directly and we are not able to give a conclusion. But the computer experiment indicates that they are probably quasi-disks.

\subsection{Organization of the paper}

To study the connectedness locus and capture components of $\MS(\theta)$, we shall consider the critically marked parameter space $\MP^{cm}(\theta)$ of cubic polynomials having a fixed Siegel disk at $0$ with rotation number $\theta$ (see \S\ref{subsec:basic}). Such a parameterization was introduced by Zakeri in \cite{Zak99}. The advantage of introducing $\MP^{cm}(\theta)$ is that one can trace which critical point lies on the boundary of the Siegel disk. The Main Theorem will be reduced to prove that all capture components in $\MP^{cm}(\theta)$ are Jordan domains.

In \S\ref{sec:parameterization}, we first recall some basic definitions and known results of $\MP^{cm}(\theta)$, and then construct a conformal isomorphism between each capture component and the unit disk by quasiconformal surgery. This allows us to define the parameter rays in capture components.

In \S\ref{sec:Jordan}, based on the continuous dependence of the boundaries of the Siegel disks, we prove that the impression of each parameter ray is a singleton by holomorphic motion and a rigidity argument. Combining the results of cubic Siegel polynomials in dynamical planes, we prove that the boundary of each capture component is a Jordan curve.

In \S\ref{sec:quasi-disk}, we prove that a capture component is a quasi-disk if the corresponding boundary is disjoint with the Zakeri curve, and the idea is essentially due to Shishikura who proved the similarity between the parameter space and dynamical planes by holomorphic motion.

In \S\ref{sec:num-cap}, we calculate the number of capture components of any given level, and it turns out that it suffices to count the number of centers in capture components.

\medskip
The method in this paper does not depend on puzzles but relies on the continuous variation of the Siegel disk boundaries. Hence one can obtain the same result for the parameter spaces of quadratic rational maps having a bounded type Siegel disk of period one \cite{Zha08b} and period two \cite{FYZ22}.

\vskip0.2cm
\noindent\textbf{Notations.}
Let $\C$ be the complex plane, $\EC=\C\cup\{\infty\}$ the Riemann sphere and denote $\C^*=\C\setminus\{0\}$.
For $a\in\C$ and $r>0$, we denote $\D(a,r)=\{z\in\C:|z-a|<r\}$ and let $\D=\D(0,1)$ be the unit disk.
For a Jordan curve $\gamma$ in $\C$, we denote by $\gamma_\Ext$ and $\gamma_\Int$ the unbounded and bounded connected components of $\C\setminus\gamma$ respectively.

\vskip0.2cm
\noindent\textbf{Acknowledgements.}
We would like to thank Arnaud Ch\'{e}ritat for sharing the idea to draw the bifurcation loci in this paper, and the referee for careful reading and detailed comments and corrections.
X. Wang was supported by National Key R\&D Program of China (Grants No. 2021YFA1003200),  NSFC (Grant No. 12131016), and the Fundamental Research Funds for the Central Universities  226-2024-00136. F. Yang was supported by NSFC (Grant No. 12571093).

\section{Parametrization of capture components} \label{sec:parameterization}

In this section we first recall some definitions and basic results in \cite{Zak99} and then give a parameterization of all capture components. In the rest of this paper, we fix the bounded type irrational number $\theta$.

\subsection{Critical marked cubic Siegel polynomials}\label{subsec:basic}

\begin{figure}[!htpb]
  \setlength{\unitlength}{1mm}
  \centering
  \includegraphics[width=0.95\textwidth]{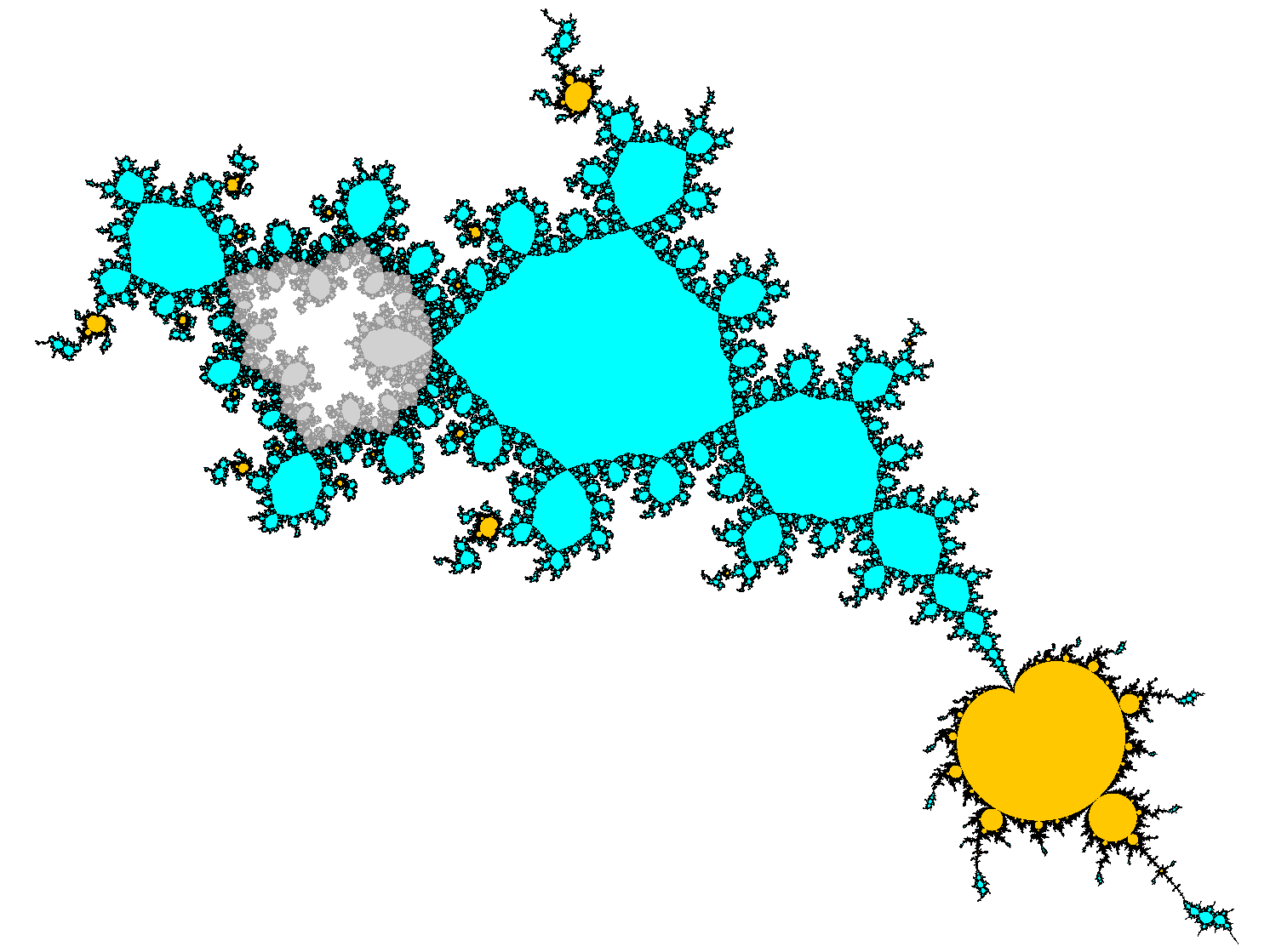}
  \caption{The parameter plane of $\MP^{cm}(\theta)$ with $\theta=(\sqrt{5}-1)/2$. The capture components outside the Zakeri curve are colored cyan and the hyperbolic-like components are colored yellow. The two parameter spaces $\MP^{cm}(\theta)$ and $\MS(\theta)$ are related by \eqref{equ:relation-a-c}. See Figure \ref{Fig:para-BH}.}
  \label{Fig:parameter}
\end{figure}

Following Zakeri, we consider the space of \textit{critically marked} cubic Siegel polynomials with rotation number $\theta$:
\begin{equation}\label{equ:crit-marked}
\mathcal{P}^{cm}(\theta):=\left\{P_c(z)=e^{2\pi\ii\theta}z\left(1-\frac{1}{2}\Big(1+\frac{1}{c}\Big)z+\frac{z^2}{3c}\right):c\in\C^*\right\}.
\end{equation}
Every $P_c\in\MP^{cm}(\theta)$ has two critical points $1$ and $c$. Let $\Delta_c=\Delta(P_c)$ be the fixed Siegel disk of $P_c$ centered at $0$ with rotation number $\theta$. It is known that $\partial\Delta_c$ passes through at least one critical point $1$ or $c$ (see \cite{Zak99} or \cite{GS03}). We denote by
\begin{equation}\label{equ:Zakeri-curve}
\Gamma:=\{c\in\C^*:\,\partial \Delta_c \text{ contains both } 1 \text{ and }c\}.
\end{equation}
The set $\Gamma$ is referred as the \textit{Zakeri curve}.
The forward orbit of $z\in\C$ under $P_c$ is denoted as $O^+(z):=\bigcup_{k\geqslant 0}P_c^{\circ k}(z)$.
The \textit{capture set} and the \textit{connectedness locus} of $\MP^{cm}(\theta)$, respectively, are defined as
\begin{equation}
\begin{split}
\mathcal{C}ap:=&~\big\{c\in\C^*:\big(O^+(1)\cup O^+(c)\big)\cap \Delta_c\neq \emptyset\big\}, \text{~~ and} \\
\MM:=&~\big\{c\in\C^*: \text{The Julia set } J(P_c) \text{ is connected}\big\}.
\end{split}
\end{equation}
The following results were proved by Zakeri in \cite{Zak99} (see Figure \ref{Fig:parameter}).

\begin{thm}\label{thm:Zakeri}
The following properties hold:
\begin{enumerate}
\item $\Gamma$ is a Jordan curve and $\MM$ is a connected compact subset of $\C^*$ which are invariant under $c\mapsto 1/c$, where $\{1,-1\}\subset\Gamma\subset\MM\subset\{c\in\C^*:\frac{1}{30}<|c|<30\}$;
\item $\partial\Delta_c$ is a quasi-circle passing through only one critical point $1$ if $c\in\Gamma_\Ext$, passing through only one critical point $c$ if $c\in\Gamma_\Int$ and passing through both if $c\in\Gamma$;
\item $\C^*\ni c\mapsto\partial\Delta_c$ is continuous in the Hausdorff topology, and moreover, $\partial\Delta_c$ moves holomorphically as $c\in\Gamma_\Ext$ and $\Gamma_\Int$; and
\item all connected components of $\mathcal{C}ap$ are open and contained in $\MM\setminus\Gamma$.
\end{enumerate}
\end{thm}

Every connected component of $\mathcal{C}ap$ is called a \textit{capture component}.
Let $\MS(\theta)=\{f_a(z)=e^{2\pi\ii\theta}z+a z^2+z^3:a\in\C\}$ be introduced in the introduction.

\begin{lem}\label{lem:covering}
If all capture components in $\MP^{cm}(\theta)$ are Jordan domains, then all capture components in $\MS(\theta)$ are also.
\end{lem}

\begin{proof}
A direct calculation shows that $P_c\in \MP^{cm}(\theta)$ is conformally conjugate to $f_a\in\MS(\theta)$ if and only if
\begin{equation}\label{equ:relation-a-c}
a^2=\eta(c):=\frac{3 e^{2\pi\ii\theta}}{4}\left(c+\frac{1}{c}+2\right).
\end{equation}
By Theorem \ref{thm:Zakeri}(a), the Zakeri curve $\Gamma\subset\C^*$ is a Jordan curve containing $\pm 1$ which is invariant under $c\mapsto 1/c$. This implies that
$c\mapsto c+\frac{1}{c}:\C^*\to\C$ maps $\Gamma$ onto a simple arc connecting $2$ with $-2$. Hence $\eta(\Gamma)$ is a simple arc connecting $3 e^{2\pi\ii\theta}$ with $0$. Thus $\eta:\Gamma_{\Ext}\to\C\setminus\eta(\Gamma)$ and $\eta:\Gamma_{\Int}\to\C\setminus\eta(\Gamma)$ are conformal.

Let $\MC$ be a capture component in $\MP^{cm}(\theta)$. By Theorem \ref{thm:Zakeri}(d), if $\MC$ is a Jordan domain, then $\eta(\MC)$ is a Jordan domain avoiding $0$. Hence $\{a\in\C:a^2\in\eta(\MC)\}$ consists of two disjoint Jordan domains which are capture components in $\MS(\theta)$.
\end{proof}

By the symmetry of $\MP^{cm}(\theta)$ as stated in Theorem \ref{thm:Zakeri}(a), we only need to consider the parameters outside the Zakeri curve $\Gamma$. Without loss of generality, in the following we assume that
\begin{equation}
c\in\Gamma\cup\Gamma_\Ext.
\end{equation}
Under this assumption, from Theorem \ref{thm:Zakeri}(b) we know that $\partial\Delta_c$ contains the critical point $1$. Moreover, the parameter $c$ is contained in a capture component in $\Gamma_\Ext$ if and only if $O^+(c)\cap\Delta_c\neq\emptyset$.

\subsection{Parameterization}

The concept of holomorphic motion was introduced by Ma\~{n}\'{e}-Sad-Sullivan and Lyubich independently (see \cite{MSS83},  \cite{Lyu83b}).

\begin{defi}[{Holomorphic motion}]
Let $E$ be a subset of $\EC$ and $\Lambda$ a connected complex manifold. A map $h:\Lambda\times E\rightarrow\EC$ is called a \textit{holomorphic motion} of $E$ parameterized by $\Lambda$ with base point $\lambda_0$ if
\begin{itemize}
\item $h(\lambda_0,z)=z$ for all $z\in E$;
\item for every $z\in E$, $\lambda\mapsto h(\lambda,z)$ is holomorphic in $\Lambda$; and
\item for every $\lambda\in \Lambda$, $z\mapsto h(\lambda,z)$ is injective on $E$.
\end{itemize}
\end{defi}

A basic fact about the holomorphic motion is the $\lambda$-\textit{lemma}: any holomorphic motion $h:\Lambda\times E\rightarrow\EC$ can be extended to a holomorphic motion $h:\Lambda\times \overline{E}\rightarrow\EC$. The following more general result can be found in \cite{Slo91}, \cite{GJW10}.

\begin{thm}\label{thm:holo-extension}
Let $h:\mathbb{D}\times E\rightarrow\EC$ a holomorphic motion of a subset $E$ of $\EC$ parameterized by $\D$ with base point $0$. Then there exists a holomorphic motion $H:\D\times\EC\to\EC$ which extends $h:\mathbb{D}\times E\rightarrow\EC$. Moreover, for any fixed $\lambda\in\D$, $H_\lambda=H(\lambda,\cdot):\EC\to\EC$ is a quasiconformal map whose complex dilatation satisfies
\begin{equation}
K(H_\lambda)\leqslant \frac{1+|\lambda|}{1-|\lambda|}.
\end{equation}
\end{thm}

For $c\in \Gamma\cup\Gamma_\Ext$, the boundary of the Siegel disk $\Delta_c$ of $P_c$ contains the critical point $1$. Let $\phi_c:\Delta_c\to\D$ be the unique conformal isomorphism whose homeomorphic extension $\phi_c:\overline{\Delta}_c\to\overline{\D}$ satisfies
\begin{equation}\label{equ:linear-map}
\phi_c(P_c(z))=e^{2\pi\ii\theta}\phi_c(z) \text{ for all } z\in\overline{\Delta}_c \text{\quad and\quad}\phi_c(1)=1.
\end{equation}
In other words, $\phi_c$ is the normalized linearizing map of $P_c$ in $\overline{\Delta}_c$ which maps the critical point $1\in\partial\Delta_c$ to $1\in\partial\D$.

If $c_0\in\Gamma_\Ext$ lies in a capture component $\MC$, then there exists a unique integer $\ell\geqslant 1$ such that
$P_{c_0}^{\circ (\ell-1)}(c_0)\not\in \Delta_{c_0}$ and $P_{c_0}^{\circ \ell}(c_0)\in \Delta_{c_0}$. By Theorem \ref{thm:Zakeri}(c), we have
\begin{equation}
P_{c}^{\circ (\ell-1)}(c)\not\in \Delta_{c} \text{ and } P_{c}^{\circ \ell}(c)\in \Delta_{c} \text{\quad for all } c\in\MC.
\end{equation}
The integer $\ell\geqslant 1$ is called the \textit{level} of the capture component $\MC$.

\begin{lem} \label{para-c}
Let $\MC$ be a capture component in $\Gamma_\Ext$ of level $\ell\geqslant 1$. Then the map
\begin{equation}
\Phi:
\begin{cases}
\MC\rightarrow \mathbb{D}\ \  &\\
c\mapsto \phi_c(P_c^{\circ \ell}(c))\ \ &
\end{cases}
\end{equation}
is a conformal isomorphism.
\end{lem}

\begin{proof}
By Theorem \ref{thm:Zakeri}(c), the boundary of the Siegel disk of $P_c$ moves holomorphically with respect to $c\in\Gamma_\Ext$. Hence there exists a holomorphic motion of $\partial\Delta_c$ parameterized by $\Gamma_\Ext$.
According to \cite[Main Theorem, (iv) $\Rightarrow$ (i), p.\,1045]{Zak16} and \cite[Lemma 6.1(2)]{Zak16}, the normalized linearizing map $\phi_c(z)$ depends holomorphically on $c$ for any fixed $z\in\Delta_c$. Hence $\Phi(c)=\phi_c(P_c^{\circ \ell}(c))$ is holomorphic in $\MC$.
To show that $\Phi$ is a conformal map, a standard method is to prove that $\Phi$ is proper and that there exists a unique $c\in\MC$ such that $\Phi(c)=0$. The former has been proved in \cite[Lemma 5.4]{Zak99} and the latter follows from \cite[Theorem 5.5]{Zak99}.
Here we use a different method: constructing a holomorphic map $\Psi: \mathbb{D}\rightarrow \MC$ such that $\Phi\circ \Psi =\id$, which implies that $\Phi:\MC\to\D$ is conformal.

The idea of the proof is to use quasiconformal surgery. Fix $c_0\in \MC$ and set $\zeta_0=\Phi(c_0)$. For $r>0$, let $D(\zeta_0,r)=\{\zeta\in\mathbb{D}: d_{\mathbb{D}}(\zeta,\zeta_0)<r\}$ be the hyperbolic disk centered at $\zeta_0$ with radius $r$. Let $E=\partial D(\zeta_0,r)\cup\{\zeta_0\}$. We define a map $h: D(\zeta_0, r)\times E\rightarrow\mathbb{D}$ by:
\begin{equation}
h(\zeta,z):=
\begin{cases}
 z,\ \  &z\in \partial D(\zeta_0, r),\\
\zeta,\ \ &z=\zeta_0.
\end{cases}
\end{equation}
It is easy to check that
\begin{itemize}
\item $h(\zeta_0,z)=z$ for all $z\in E $;
\item for every fixed $\zeta\in D(\zeta_0, r), \ z\mapsto h(\zeta,z)$ is injective on $E$; and
\item for every fixed $z\in E, \ \zeta\mapsto h(\zeta,z)$ is holomorphic in $D(\zeta_0, r)$.
\end{itemize}
Thus $h: D(\zeta_0, r)\times E\rightarrow\mathbb{D}$ is a holomorphic motion parameterized by $D(\zeta_0, r)$ with base point $\zeta_0$.

By Theorem \ref{thm:holo-extension}, $h$ admits an extension $H: D(\zeta_0, r)\times \mathbb{\widehat{C}}\rightarrow\mathbb{\widehat{C}}$, which is also a holomorphic motion.
Let $B$ be the Fatou component of $P_{c_0}$ containing $P_{c_0}^{\circ(\ell-1)}(c_0)$ and $B_{r}=(P_{c_0}|_{B})^{-1}\circ(\phi_{c_0})^{-1}(D(\zeta_0,r))\subset B$.
Then $P_{c_0}^{\circ(\ell-1)}(c_0)\in B_r$ since $\zeta_0=\phi_{c_0}(P_{c_0}^{\circ\ell}(c_0))$. We define
\begin{equation}
\psi_{\zeta}(z):=\phi_{c_0}^{-1}\circ H(\zeta,\phi_{c_0}\circ P_{c_0}(z)),
\end{equation}
where $(\zeta,z)\in D(\zeta_0, r)\times B_r$. It is easy to check that $\psi_{\zeta}$ satisfies:
\begin{itemize}
\item $\psi_{\zeta_0}(z)=P_{c_0}(z)$ for all $z\in B_{r}$;
\item $\psi_\zeta(z)=P_{c_0}(z)$ for any fixed $\zeta$ and for all $z\in\partial B_r$;
\item $\psi_{\zeta}(P_{c_0}^{\circ(\ell-1)}(c_0))=\phi_{c_0}^{-1}(\zeta)$ for all $\zeta\in D(\zeta_0, r)$;
\item $\psi_\zeta: B_{r}\rightarrow (\phi_{c_0})^{-1}(D(\zeta_0,r))$ is a quasi-regular\footnote{It is quasiconformal if the level $\ell\geqslant 2$ since in this case $B_r$ does not contain the critical point $c_0$.} map for any fixed $\zeta$; and
\item $\zeta\mapsto\psi_{\zeta}(z)$ is holomorphic in $D(\zeta_0, r)$ for any fixed $z\in B_{r}$.
\end{itemize}

Now we define a quasi-regular map
\begin{equation}
Q_\zeta(z):=\begin{cases}
 \psi_{\zeta}(z),\ \  & z\in B_{r},\\
P_{c_0}(z),\ \ &z\in \widehat{\mathbb{C}}\setminus B_{r}.
 \end{cases}
 \end{equation}
Let $\sigma$ be the standard complex structure. We construct a $Q_\zeta$-invariant complex structure $\sigma_\zeta$ as follows:
\begin{equation}
\sigma_\zeta:=
\begin{cases} (Q_\zeta^{\circ n})^*(\psi_{\zeta}^*\sigma),\ \
 &\mathrm{in} \ P_{c_0}^{-n}(B_r), \  \ n\geqslant 0,\\
 \sigma,\ \ &\mathrm{in}\ \widehat{\mathbb{C}}\setminus \bigcup_{n\geqslant0} P_{c_0}^{-n}(B_{r}).
 \end{cases}
 \end{equation}

The Beltrami coefficient $\mu_\zeta$ of  $\sigma_\zeta$ satisfies $\|\mu_\zeta\|<1$ for all $\zeta\in D(\zeta_0, r)$ since $Q_\zeta$ is holomorphic outside $B_r$.
By the Measurable Riemann Mapping Theorem, there exists a family of quasiconformal maps $\varphi_\zeta:\mathbb{\widehat{C}}\rightarrow \mathbb{\widehat{C}}$, parameterized by $\zeta\in D(\zeta_0,r)$, each  solves the Beltrami equation $\overline{\partial} \varphi_\zeta= \mu_\zeta {\partial} \varphi_\zeta$, and fixes $0,1,\infty$.
Then there is a holomorphic map $\Psi: D(\zeta_0, r)\rightarrow\MC$ such that $\Psi(\zeta_0)=c_0$ and $\varphi_\zeta\circ Q_\zeta\circ \varphi_\zeta^{-1}=P_{\Psi(\zeta)}$.

Note that $\Psi(\zeta)$ is a critical point of $P_{\Psi(\zeta)}$ for every $\zeta\in D(\zeta_0, r)$. Hence $\varphi_\zeta^{-1}(\Psi(\zeta))=c_0$.
Since both $\phi_{c_0}\circ\varphi_\zeta^{-1}$ and $\phi_{\Psi(\zeta)}$ are normalized linearizing maps of $P_{\Psi(\zeta)}$ in $\overline{\Delta}_{\Psi(\zeta)}$ satisfying \eqref{equ:linear-map}, we have $\phi_{c_0}\circ\varphi_\zeta^{-1}=\phi_{\Psi(\zeta)}$.
Therefore, for $\zeta\in D(\zeta_0, r)$, we have
\begin{equation}
\begin{split}
\Phi(\Psi(\zeta))=~&
\phi_{\Psi(\zeta)}(P_{\Psi(\zeta)}^{\circ \ell}(\Psi(\zeta)))
=\phi_{\Psi(\zeta)}\circ \varphi_\zeta\circ Q_\zeta^{\circ\ell} \circ \varphi_\zeta^{-1}(\Psi(\zeta))\\
=~&\phi_{\Psi(\zeta)}\circ \varphi_\zeta\circ \psi_\zeta\circ
P_{c_0}^{\circ (\ell-1)}(c_0)=\phi_{\Psi(\zeta)}\circ \varphi_\zeta\circ \phi_{c_0}^{-1}(\zeta)\\
=~&\phi_{\Psi(\zeta)}\circ \phi_{\Psi(\zeta)}^{-1}(\zeta)=\zeta.
\end{split}
\end{equation}
Since $r$ can be arbitrarily large,  the map $\Phi$ actually admits a global inverse.
\end{proof}

\begin{rmk}
Each capture component $\MC$ (with level $\ell$) contains exactly one \textit{center}, i.e., there exists a unique $c\in\MC$ such that $P_c^{\circ \ell}(c)=0$.
\end{rmk}

\section{Capture components are Jordan domains}\label{sec:Jordan}

For $c\in\Gamma\cup\Gamma_\Ext$, let $U_c$ be the Fatou component (if any) of $P_c$ containing the critical point $c$. Let $U$ be any connected component of $\bigcup_{n\geqslant 0}P_c^{-n}(\Delta_c)\setminus \bigcup_{k\geqslant 0}P_c^{-k}(U_c)$. Suppose that $n\geqslant 0$ is the smallest integer such that $P_c^{\circ n}(U)=\Delta_c$. The \textit{dynamical internal ray} $R_U(t)$ of angle $t\in \mathbb{R}/\mathbb{Z}$ in $U$ is defined as
\begin{equation}
R_U(t):=(P_c^{\circ n}|_U)^{-1}(\phi_c^{-1}((0,1)e^{2\pi\ii t})).
\end{equation}
According to Theorem \ref{thm:Zakeri}(b), the Siegel disk $\Delta_c$ and its all preimages are Jordan domains. This means that all dynamical internal rays land.

\medskip
Let $\MC\subset\Gamma_\Ext$ be a capture component with level $\ell\geqslant 1$. For $c\in\MC$, we denote
\begin{equation}
V_c:= \text{The Fatou component of } P_c \text{ containing } v_c,
\end{equation}
where $v_c=P_c(c)$.
Clearly, the {\it center} $\xi=\xi(c)$ of  $V_{c}$, defined as  the unique point $\xi\in  V_{c}$ satisfying $P_{c}^{\circ (\ell-1)}(\xi)=0$,   moves continuously with respect to   $c\in\MC$.  The center map $c\mapsto \xi(c)$ has a continuous extension to $\partial\MC$. Therefore,  when $c\in \partial\MC$, the point $\xi(c)$ and the Fatou component  $V_{c}$ containing $\xi(c)$ are  well-defined. In particular, when $c\in \partial\MC$, we still have $P_c^{\circ (\ell-1)}(V_c)=\Delta_c$ and $P_c^{\circ (\ell-2)}(V_c)\neq\Delta_c$ (if $\ell\geqslant 2$).

\medskip
By Lemma \ref{para-c}, the map
 $$\Phi: \MC\rightarrow \mathbb{D}, \  {c}\mapsto \phi_{c}(P_{c}^{\circ \ell}(c))$$
is a conformal isomorphism.
The {\it parameter ray} $\mathcal{R}(t)$ of angle $t\in \mathbb{R}/\mathbb{Z}$ in $\MC$ is
\begin{equation}
\mathcal{R}(t):=\Phi^{-1}((0,1)e^{2\pi \ii t}).
\end{equation}
The  {\it impression} $\mathcal{I}(t)$ of $\mathcal{R}(t)$ is defined as $\mathcal{I}(t)=\bigcap_{k\geqslant 1}\overline{\mathcal{S}_k(t)}$, where
\begin{equation}
\mathcal{S}_k(t):=\Phi^{-1}\Big(\big\{re^{2\pi \ii \alpha}:\, r\in \big(1-\tfrac{1}{k},1\big), \ \alpha\in\big(t-\tfrac{1}{k},t+\tfrac{1}{k}\big)\big\}\Big).
\end{equation}

 \begin{lem}\label{dyn-para}
 For any $t\in  \mathbb{R}/\mathbb{Z}$ and any $c_0\in \mathcal{I}(t)$, the dynamical internal ray $R_{V_{c_0}}(t)$ lands at $v_{c_0}\in\partial V_{c_0}$.
 \end{lem}

\begin{proof}
By Theorem \ref{thm:Zakeri}(c), the boundary $\partial \Delta_c$ moves continuously in the Hausdorff topology with respect to $c\in \mathbb{C}^*$.  It follows that
$P_c^{-(\ell-1)}(\partial \Delta_c)$ moves continuously in the Hausdorff topology with respect to $c\in \mathbb{C}^*$.
Therefore, the maps $c\mapsto v_c$ and $c\mapsto\partial V_c$ are continuous with respect to $c\in\overline{\MC}$.
Note that if $c=\Phi^{-1}(re^{2\pi\ii t})\in \mathcal{R}(t)$ with $r\in(0,1)$, then
\begin{equation}
v_c=P_c(c)=(P_c^{\circ (\ell-1)}|_{V_c})^{-1}(\phi_c^{-1}(re^{2\pi\ii t}))\in R_{V_c}(t).
\end{equation}
Hence if $c_0\in \mathcal{I}(t)$, then $v_{c_0}\in\partial V_{c_0}$. Thus if $c\mapsto R_{V_c}(t)$ moves continuously with respect to $c\in\overline{\MC}$, then the dynamical internal ray $R_{V_{c_0}}(t)$ lands at $v_{c_0}$.

To prove the continuity of $c\mapsto R_{V_c}(t)$, it suffices to prove that $c\mapsto \phi_c^{-1}(re^{2\pi \ii t})$ moves continuously with respect to $c\in\overline{\MC}$ for given $r\in(0,1)$ and $t\in\mathbb{R}/\mathbb{Z}$.
In fact, by \eqref{equ:linear-map}, we have
\begin{equation}
P_c(\phi_c^{-1}(z))=\phi_c^{-1}(e^{2\pi\ii\theta}z) \text{ for all } z\in\overline{\D} \text{\quad and\quad}\phi_c^{-1}(1)=1.
\end{equation}
By comparing the coefficients, the conformal map $\phi_c^{-1}:\D\to\Delta_c$ can be written as the Taylor expansion $\phi_c^{-1}(z)=\sum_{n=1}^\infty a_n(c) z^n:\D\to\Delta_c$, where $c\mapsto a_n(c)$ is continuous with respect to $c\in\overline{\MC}$ for all $n\geqslant 1$.
Since $c\mapsto\partial\Delta_c$ moves continuously, it follows that there exists $M_0\geqslant 1$ such that $|a_1(c)|\leqslant M_0$ for all $c\in\overline{\MC}$. Since $\phi_c^{-1}:\D\to\Delta_c$ is conformal, by \cite[Theorem 1.8, p.\,4]{CG93} we have
\begin{equation}
|a_n(c)|\leqslant |a_1(c)| en^2\leqslant M_0 en^2 \text{\quad for all } n\geqslant 1 \text{ and } c\in\overline{\MC}.
\end{equation}
Let $\varepsilon>0$, $c_0\in\overline{\MC}$ and $z_0=re^{2\pi\ii t}\in\D\setminus\{0\}$ be given. There exists $N\geqslant 1$ such that if $c\in\overline{\MC}$ is sufficiently close to $c_0$, then
\begin{equation}
\big|\phi_c^{-1}(z_0)-\phi_{c_0}^{-1}(z_0)\big|\leqslant\sum_{n=1}^N |a_n(c)-a_n(c_0)|r^n+2\sum_{n=N+1}^\infty M_0 en^2 r^n<\frac{\varepsilon}{3}+\frac{2\varepsilon}{3}=\varepsilon.
\end{equation}
Therefore, for given $r\in(0,1)$ and $t\in\mathbb{R}/\mathbb{Z}$, the map $c\mapsto \phi_c^{-1}(re^{2\pi \ii t})$ is continuous with respect to $c\in\overline{\MC}$ and the proof is complete.
\end{proof}

Recall that $\Gamma=\{c\in\C^*:\,1,c\in\partial\Delta_c\}$ is the Zakeri curve.

\begin{lem}\label{crit-relation}
For any $t\in  \mathbb{R}/\mathbb{Z}$ such that $\mathcal{I}(t)\cap\Gamma\neq \emptyset$, we have $t=\theta$.
Further, for any $c\in \mathcal{I}(t)\cap\Gamma$, we have $P_c^{\circ (\ell-1)}(c)=1$.
\end{lem}

\begin{proof}
Let $c\in \mathcal{I}(t)\cap\Gamma$. Then we have $c\in \partial \Delta_c$, which implies that $P_c^{\circ k}(c)\in \partial \Delta_c$ for all
$k\geqslant 0$. We first assume $\ell\geqslant 2$. By Lemma \ref{dyn-para}, we have
$P_c^{\circ (\ell-1)}(c)\in \partial P_c^{\circ (\ell-2)}(V_c)$.
Therefore, $P_c^{\circ (\ell-1)}(c)$ lies on the common boundary of $P_c^{\circ (\ell-2)}(V_c)$ and $\Delta_c$.
The fact $P_c(P_c^{\circ (\ell-2)}(V_c))=\Delta_c=P_c(\Delta_c)$ implies that $P_c^{\circ (\ell-1)}(c)$ is a critical point.
This critical point cannot be $c$ since otherwise $c$ would be a super-attracting period point, which contradicts $c\in\partial\Delta_c$. So we have
$P_c^{\circ (\ell-1)}(c)=1$. By the fact that $P_c(R_{\Delta_c}(s))=R_{\Delta_c}(s+\theta)$ for all $s\in\R/\Z$ and the normalization of $\phi_c$ in \eqref{equ:linear-map}, we conclude that $P_c^{\circ \ell}(c)$ is the landing point of the dynamical internal ray $R_{\Delta_c}(\theta)$. Hence $t=\theta$.

Suppose $\ell=1$. In this case we have $V_c=\Delta_c$ for all $c\in\overline{\MC}$. For $c\in\MC$, recall that $U_c$ is the Fatou component of $P_c$ containing $c$. Note that $P_c:U_c\to V_c=\Delta_c$ is a branched covering of degree two. Since $\overline{\Delta}_c$ moves continuously with respect to $c\in\overline{\MC}$, if $c\in\MI(t)\subset\partial\MC$, there exist two Fatou components $U_c'$ and $U_c''$ which are different from $\Delta_c$ such that $c\in\partial U'\cap\partial U_c''$ and $P_c(U_c')=P_c(U_c'')=\Delta_c$. This implies that $c$ is a double critical point if further $c\in\Gamma$. Hence $c=1$ and $P_c(c)=P_c(1)$ is the landing point of the internal ray $R_{\Delta_c}(\theta)$. Hence we also have $t=\theta$.
\end{proof}

Let $\MM$ be the connectedness locus of $\MP^{cm}(\theta)$. For each $c\in\MM$, we use $\Omega_c$ to denote the basin of infinity of $P_c$ and denote by
\begin{equation}
\chi_c:\Omega_c\to\EC\setminus\overline{\D}
\end{equation}
the unique B\"{o}ttcher map satisfying $\chi_c(P_c(z))=(\chi_c(z))^3$ and $\chi_c'(\infty)=1$.

\begin{lem}\label{lc}
For any $t\in  \mathbb{R}/\mathbb{Z}$, the impression $\mathcal{I}(t)$ is a singleton.
\end{lem}

\begin{proof}
Let $\mathcal{I}^*(t)= \mathcal{I}(t)\setminus\Gamma$. We will first show that $\mathcal{I}^*(t)$ is a singleton or empty.
If it is not true, then there exists a connected compact subset $\mathcal{E}$ of $\mathcal{I}^*(t)$  containing at least two points.
By Lemma  \ref{dyn-para}, the dynamical internal ray $R_{V_{c}}(t)$ lands at $v_{c}$ for all $c\in \mathcal{E}$.
By the continuity in Theorem \ref{thm:Zakeri}(c) and shrinking $\mathcal{E}$ if necessary, we may assume that all $c\in\mathcal{E}$ satisfy $P_{c}^{\circ (\ell-1)}(c)\notin \overline{\Delta}_c$ and $P_{c}^{\circ \ell}(c)\in \partial \Delta_{c}$.
Then there is a Jordan domain $\mathcal{D}$ with $\mathcal{E}\subset \mathcal{D}\subset \mathbb{C}^*\setminus\Gamma$, so that for all $c\in \mathcal{D}$, we have $P_{c}^{\circ (\ell-1)}(c)\notin \overline{\Delta}_{c}$.

For $c_1\in\mathcal{E}$, we denote
\begin{equation}
W_{c_1}:=\{P_{c_1}^{\circ j}(c_1):\, 1\leqslant j \leqslant \ell-1\}\cup \overline{\Delta}_{c_1}\cup\Omega_{c_1}.
\end{equation}
It is clear that ${W_{c_1}}$ contains the post-critical set of $P_{c_1}$.
We define a continuous map $h: \mathcal{D}\times W_{c_1}\rightarrow \mathbb{\widehat{C}}$ satisfying that
\begin{itemize}
\item $h(c, P_{c_1}^{\circ j}({c_1}))=P_{c}^{\circ j}(c) $ for all $c\in \mathcal{D}$ and  $1\leqslant j\leqslant \ell-1$;
\item $h(c, z)=\phi_{c}^{-1}\circ \phi_{c_1}(z)$ for all $(c,z)\in \mathcal{D}\times \overline{\Delta}_{c_1}$; and
\item $h(c, z)=\chi_{c}^{-1}\circ \chi_{c_1}(z)$ for all $(c,z)\in \mathcal{D}\times \Omega_{c_1}$.
\end{itemize}

One may verify that $h$ is a holomorphic motion, parameterized by $\mathcal{D}$ and with base point $c_1$ (i.e. $h(c_1, \cdot)=\id$). By Theorem \ref{thm:holo-extension}, there is a holomorphic motion $H:\mathcal{D}\times
\mathbb{\widehat{C}}\rightarrow\mathbb{\widehat{C}}$ extending $h$. We consider the restriction $H_0=H|_{\mathcal{E}\times\mathbb{\widehat{C}}}$ of $H$.
Note that fix any  $c\in \mathcal{E}$, the map  $z\mapsto H(c,z)$ sends the post-critical set of $P_{c_1}$ to that of $P_{c}$, preserving the dynamics on this set.
By the lifting property, there is a unique continuous map  $H_1: \mathcal{E}\times \mathbb{\widehat{C}}\rightarrow \mathbb{\widehat{C}}$ such that
$P_{c}(H_1(c,z))=H_0(c, P_{c_1}(z))$ for all $(c,z)\in\mathcal{E}\times\widehat{\mathbb{C}}$
and $H_1(c_1,\cdot)\equiv \id$. It is not hard to see that $H_{1}$ is also a holomorphic motion.

For any parameter $c_2\in \mathcal{E}$ which is different from $c_1$,  we set $\psi_0=H_0(c_2, \cdot)$ and $\psi_1=H_1(c_2, \cdot)$.
Both $\psi_0$ and $\psi_1$ are quasiconformal maps on $\EC$, satisfying  $P_{c_2}\circ \psi_{1}=\psi_0\circ P_{c_1}$.
One may verify that $\psi_0$ and $\psi_1$ are isotopic rel ${W_{c_1}}$.
Again by the lifting property, there is a sequence of quasiconformal maps $\psi_k$'s satisfying that
\begin{itemize}
\item $P_{c_2}\circ \psi_{k+1}=\psi_k\circ P_{c_1}$ for all $k\geqslant0$; and
\item $\psi_{k+1}$ and $\psi_k$ are isotopic rel $P_{c_1}^{-k}({W_{c_1}})$.
\end{itemize}

The maps $\psi_k$'s form a normal family because their dilatations are uniformly bounded above.
Any limit map  $\psi_\infty$ of $\psi_k$'s is quasiconformal in $\EC$,  holomorphic in the Fatou set $F(P_{c_1})$ of $P_{c_1}$ and  satisfies
$P_{c_2}\circ \psi_{\infty}=\psi_\infty\circ P_{c_1}$ in $F(P_{c_1})$. By the continuity, $\psi_\infty$ is a conjugacy on $\EC$.
Such $\psi_\infty$ is necessarily conformal. Indeed, otherwise $P_{c_1}$ will carries an invariant line field, hence admits a quasiconformal deformation in a neighborhood of $c_1$. This implies that $P_{c_1}$ is disjoint with the bifurcation locus and $J$-stable, which contradicts the fact that $c_1\in\partial \MC$.
Since $\psi_\infty:\EC\to\EC$ is a conformal map fixing $0$, $1$ and $\infty$ (note that $1$ is the critical point on $\partial\Delta_{c_1}$ and $\partial\Delta_{c_2}$), it implies that $\psi_\infty$ is the identity and hence $c_1=c_2$. This contradicts our assumption. Therefore, $\mathcal{I}^*(t)$ is either a singleton or empty.

\medskip
If $\MI(t)\cap\Gamma=\emptyset$, we conclude that $\MI(t)=\MI^*(t)$ is a singleton. If $\MI(t)\cap\Gamma\neq\emptyset$, by Lemma \ref{crit-relation}, we see that $t=\theta$ and
$$\mathcal{I}(\theta)\cap\Gamma\subset \{c\in \mathbb{C}^*:\, P_c^{\circ (\ell-1)}(c)=1\}.$$
Note that $P_c^{\circ (\ell-1)}(c)$ is a polynomial with respect to $c$. Hence the right side of above inclusion is a finite set.
 Since $\mathcal{I}(\theta)=\MI^*(\theta)\cup (\mathcal{I}(\theta)\cap\Gamma)$ is connected, we conclude that $\MI^*(\theta)=\emptyset$ and $\mathcal{I}(\theta)$ is a singleton.
\end{proof}

\begin{thm} \label{boundary-C}
Capture components are Jordan domains.
\end{thm}

\begin{proof}
Let $\MC$ be a capture component.
By Lemma \ref{lc}, we see that $\partial \MC$ is locally connected.
If two parameter rays $\mathcal{R}(t_1)$ and $\mathcal{R}(t_2)$ land at the same point $c\in \partial \MC$, then by Lemma \ref{dyn-para}, the dynamical internal rays   $R_{V_{c}}(t_1)$ and $R_{V_{c}}(t_2)$ both land at the critical value $v_{c}\in\partial V_{c}$. Since  $\partial V_c$  is a Jordan curve, we have $t_1=t_2$.
This implies that $\partial\MC$ is a Jordan curve.
\end{proof}

\begin{proof}[Proof of the Main Theorem]
This follows from Theorem \ref{boundary-C} and Lemma \ref{lem:covering}.
\end{proof}

From the proofs of Lemmas \ref{crit-relation} and \ref{lc}, we have the following immediate result.

\begin{cor}\label{cor:singleton}
If $\partial\MC\cap \Gamma\neq \emptyset$, then $\partial \MC\cap \Gamma$ is a singleton $\{c\}$ which satisfies
\begin{enumerate}
\item $c$ is the landing point of the parameter ray $\mathcal{R}(\theta)$ in $\MC$.
\item If $\ell=1$, then $c=1$. If $\ell\geqslant 2$, then $P_c^{\circ (\ell-1)}(c)=1$ and $P_c^{\circ (\ell-2)}(c)\not\in\overline{\Delta}_c$.
\end{enumerate}
If $\partial\MC\cap \Gamma= \emptyset$, then for all $c\in\partial\MC$, $P_c^{\circ \ell}(c)\in\partial{\Delta}_c$ and $P_c^{\circ (\ell-1)}(c)\not\in\overline{\Delta}_c$.
\end{cor}

\section{Some capture components are quasi-disks}\label{sec:quasi-disk}

In this section we prove that a capture component $\MC$ is a quasi-disk if $\partial \MC$ is disjoint with the Zakeri curve $\Gamma$.
To this end, we need the following Lemma \ref{lem:holo-motion}, which establishes a ``similarity" between the parameter space and the dynamical planes. The proof is strongly inspired by \cite[Lemma 3.2]{Shi98} (see also \cite[Remark 3.1]{Shi98}). For completeness we include a proof here.

For $z_0\in\C$ and $r>0$, let $\D(z_0,r)=\{z\in\C:|z-z_0|<r\}$ be the Euclidean disk of radius $r$ centered at $z_0$ and denote $\D_r:=\D(0,r)$ for short.
Let $B_r(z_0)$ be the disk of radius $r$ centered at $z_0$ in the spherical metric.

\begin{lem}\label{lem:holo-motion}
Let $\MD=\D(\lambda_0,\widetilde{s}_0)\subset\C$, $X$ be a subset of $\EC$ and $h:\MD\times X\to\EC$ a holomorphic motion parameterized by $\MD$ with base point $\lambda_0$. Suppose $v:\MD\to\EC$ is an analytic map satisfying $v(\lambda_0)=z_0\in X$ and $v(\cdot)\not\equiv h(\cdot,z_0)$. Then there exist two constants $r_0>0$ and $0<s_0<\widetilde{s}_0$ such that
\begin{enumerate}
\item For all $0<r\leqslant r_0$,
\begin{equation}
Y^r:=\{\lambda\in\D(\lambda_0,s_0):\,v(\lambda)\in h_\lambda(X\cap B_r(z_0))\}
\end{equation}
 is mapped onto $X\cap B_r(z_0)$ by the restriction of a quasi-regular map\footnote{The domain of definition of $\varphi_r=v\circ G_r^{-1}$ is seen as $G_r(\D(\lambda_0,\widetilde{s}_0))$.} $\varphi_r=v\circ G_r^{-1}$, where $G_r:\EC\to\EC$ is a $K_r$-quasiconformal mapping;
\item $K_r\to 1$ as $r\to 0$ and $(\varphi_r)_{r\in(0,r_0]}$ are quasiconformal if $v'(\lambda_0)\neq 0$.
\end{enumerate}
\end{lem}

\begin{proof}
Without loss of generality, we assume that $X$ contains at least two points, $\MD=\D$, $z_0=0$ and $h_\lambda(0)\equiv 0$ since one can change the coordinates by M\"{o}bius transformations depending analytically.

\vskip0.1cm
\textbf{Case 1.} We first assume that $v'(0)\neq 0$. Then there exist two positive constants $a_0$ and $\rho$ such that $v$ is univalent in $\D_\rho$ and
\begin{equation}
a_0|\lambda|\leqslant |v(\lambda)|<\infty, \text{\quad for all } \lambda\in\overline{\D}_\rho.
\end{equation}
By Theorem \ref{thm:holo-extension}, $h:\D\times X\to \EC$ can be extended to a holomorphic motion $H:\D\times \EC\to \EC$. For every $r>0$, we denote by
\begin{equation}
b_r:=\sup\{|H_\lambda(z)|:\,z\in \D_r \text{ and }|\lambda|\leqslant \rho\}.
\end{equation}
Since $H$ is continuous, we have $b_r\to 0$ as $r\to 0$. Then there exists a constant $r_0>0$ such that $a_0\rho\geqslant 2b_r$ for all $0<r\leqslant r_0$ and that $v(\D_\rho)\supset\overline{\D}_{r_0}$.

\vskip0.1cm
We take an $r\in(0,r_0]$. For any given $z\in X\cap\D_r$ and $\mu\in\D_{R_r}$ with $R_r:=a_0\rho/b_r\geqslant 2$, we consider the following equation with variable $\lambda$:
\begin{equation}\label{equ:holo-mu}
v(\lambda)-h_{\mu\lambda}(z)=0,
\end{equation}
where $\lambda\in D^\mu:=\big\{\lambda\in\C:|\lambda|<\min\{\rho,\tfrac{\rho}{|\mu|}\}\big\}$. Note that both $\lambda\mapsto v(\lambda)$ and $\lambda\mapsto h_{\mu\lambda}(z)$ are holomorphic in $ D^\mu$, and for $\lambda\in\partial D^\mu$, we have
\begin{equation}
|v(\lambda)|\geqslant a_0 |\lambda|=a_0\min\{\rho,\tfrac{\rho}{|\mu|}\}>b_r\geqslant |h_{\mu\lambda}(z)|.
\end{equation}
Since $v(\lambda)=0$ has a unique solution $\lambda=0$ in the disk $ D^\mu$, it follows from Rouch\'{e}'s theorem that the Equation \eqref{equ:holo-mu} has exactly one solution $\lambda=u(\mu,z)$ in $D^\mu$, and it depends analytically on $\mu$. Moreover, for $z_1$, $z_2\in X\cap\D_r$ with $z_1\neq z_2$, we have $u(\mu,z_1)\neq u(\mu,z_2)$ for all $\mu\in\D_{R_r}$ because of the injectivity of $h_\lambda$.

\vskip0.1cm
For $|\mu|<R_r$ we define
\begin{equation}
Y_\mu^r:=\{\lambda\in D^\mu:\,v(\lambda)=h_{\mu\lambda}(z) \text{ for some }z\in X\cap\D_r\}.
\end{equation}
In particular, $Y_0^r=\{\lambda\in\D_\rho:v(\lambda)\in X\cap \D_r\}$ and each point in $Y_0^r$ can be written as $v^{-1}(z)\in \D_\rho$ for $z\in X\cap\D_r$ since $v(\D_\rho)\supset\overline{\D}_r$.
Then one may verify easily that the following map is a holomorphic motion:
\begin{equation}
g^r:\D_{R_r}\times Y_0^r\to \C, \quad (\mu,v^{-1}(z))\mapsto u(\mu,z),
\end{equation}
where $z\in X\cap\D_r$. Moreover, we have $g^r(1,Y_0^r)=Y_1^r$, where
\begin{equation}
Y_1^r=\{\lambda\in\D_\rho:\,v(\lambda)\in h_\lambda(X\cap\D_r)\}.
\end{equation}
By Theorem \ref{thm:holo-extension}, $g^r:\D_{R_r}\times Y_0^r\to \C$ can be extended to a holomorphic motion $\widetilde{g}^r:\D_{R_r}\times\EC\to\EC$ so that each $\widetilde{g}_\mu^r=\widetilde{g}^r(\mu,\cdot)$ is a $K_\mu^r$-quasiconformal mapping, where $K_\mu^r=(R_r+|\mu|)/(R_r-|\mu|)$. In particular, $\widetilde{g}_1^r:\EC\to\EC$ is a $K_1^r$-quasiconformal mapping such that $\widetilde{g}_1^r(Y_0^r)=Y_1^r$, and we have $K_1^r\to 1$ as $r\to 0$ since $R_r\to\infty$ as $r\to 0$. Moreover, we have $K_1^r\leqslant 3$ for all $0<r\leqslant r_0$ since $R_r\geqslant 2$.

Since $v:Y_0^r\to X\cap\D_r$ is the restriction of the univalent map $v:\D\to\EC$, it follows that $X\cap\D_r$ is quasiconformally homeomorphic to $Y_1^r$. Hence if $v'(0)\neq 0$ then the lemma follows if we set $G_r:=\widetilde{g}_1^r$.

\vskip0.1cm
\textbf{Case 2.} Suppose $v'(0)=0$ and $v$ is of local degree $m\geqslant 2$. Note that $v\not\equiv 0$. Up to changing the coordinates, we assume further that $\infty\in X$ and $h_\lambda(\infty)=\infty$. Denote by $\omega(z)=z^m$, $\widetilde{X}_\lambda=\omega^{-1}(X_\lambda)$ and $\widetilde{X}=\omega^{-1}(X)$, where $X_\lambda=h_\lambda(X)$. By lifting $v$ and $h_\lambda$ under $\omega$, we obtain a holomorphic map $\widetilde{v}:\D\to\EC$ satisfying $v=\omega\circ \widetilde{v}$, $\widetilde{v}'(0)\neq 0$, and a holomorphic motion $\widetilde{h}_\lambda: \widetilde{X}\to \widetilde{X}_\lambda$ satisfying $h_\lambda\circ\omega=\omega\circ \widetilde{h}_\lambda$. By the arguments in Case 1, there exist $\widetilde{r}_0>0$ and $\widetilde{\rho}>0$ such that for all $0<r\leqslant \widetilde{r}_0$,
\begin{equation}
\{\lambda\in\D_{\widetilde{\rho}}:\,v(\lambda)\in h_\lambda(X\cap \D_r)\}
=\{\lambda\in\D_{\widetilde{\rho}}:\,\widetilde{v}(\lambda)\in \widetilde{h}_\lambda(\widetilde{X}\cap \omega^{-1}(\D_r))\}
\end{equation}
is mapped onto $\widetilde{X}\cap \omega^{-1}(\D_r)$ by the restriction of a quasiconformal mapping $\widetilde{\varphi}_r=\widetilde{v}\circ (\widetilde{g}_1^r)^{-1}$, where $\widetilde{v}:\D_{\widetilde{\rho}}\to\C$ is univalent, $\widetilde{g}_1^r:\EC\to\EC$ is a $K_r$-quasiconformal mapping and $K_r\to 1$ as $r\to 0$. Hence if $v'(0)= 0$ then the lemma follows if we set $G_r:=\widetilde{g}_1^r$ and $\varphi_r:=\omega\circ\widetilde{\varphi}_r=v\circ G_r^{-1}$.
\end{proof}

\begin{rmk}
The condition $v'(\lambda_0)\neq 0$, usually referred as \textit{transversality condition}, implies that there exists a similarity between the $\lambda$-parameter plane (in particular the bifurcation locus near $\lambda_0$) and the dynamical plane (in particular the Julia set near $z_0=v(\lambda_0)$). If $v'(\lambda_0)= 0$, we still have some kind of similarity and this similarity is obtained by local covering property.
\end{rmk}

For bounded type $\theta$, let $\Gamma$ be the Zakeri curve in $\MP^{cm}(\theta)$ and $\MC$ be a capture component in $\Gamma_\Ext$ with level $\ell\geqslant 1$.

\begin{thm}\label{thm:quasi-disjoint}
For any $c\in\partial\MC\setminus\Gamma$, there exists an open neighborhood $\MU$ of $c$ in $\Gamma_\Ext$ such that $\MU\cap\partial\MC$ is a quasi-arc.
In particular, $\MC$ is a quasi-disk if $\partial{\MC}\cap\Gamma=\emptyset$.
\end{thm}

\begin{proof}
Let $c_0\in\partial\MC\setminus\Gamma$, where $\MC$ is a capture component in $\Gamma_\Ext$ of level $\ell\geqslant 1$. Let $\MD=\D(c_0,\widetilde{s}_0)$ be a small neighborhood of $c_0$ such that $\overline{\MD}\cap\Gamma=\emptyset$. Consider $v(c):=P_c^{\circ\ell}(c)$. A direct calculation shows that $v(c)$ is a polynomial in $c$ (see the proof of Theorem \ref{thm:num-cap}). By Corollary \ref{cor:singleton} we have $z_0=v(c_0)\in\partial\Delta_{c_0}$. By Theorem \ref{thm:Zakeri}(c), we have a holomorphic motion $h:\MD\times\partial\Delta_{c_0}\to\C$ parameterized in $\MD$ with base point $c_0$.
Denote $h_c(\cdot):=h(c,\cdot)$. By Lemma \ref{lem:holo-motion}, there exist two positive constants $r_0$ and $s_0<\widetilde{s}_0$ such that for all $0<r\leqslant r_0$,
\begin{equation}
Y_1^r:=\{\lambda\in\D(c_0,s_0):\,v(c)\in h_c(\partial\Delta_{c_0}\cap\D_r(z_0))\}
\end{equation}
is mapped onto $\partial\Delta_{c_0}\cap \D_r(z_0)$ by the restriction of a quasi-regular map $\varphi_r=v\circ G_r^{-1}$, where $G_r:\EC\to\EC$ is  quasiconformal.

Since $\partial\Delta_{c_0}$ is a quasi-circle (see Theorem \ref{thm:Zakeri}(b)), it follows from Theorem \ref{boundary-C} that $\Comp_{c_0}(Y_1^r\cap\partial\MC)$ is a quasi-arc, where $\Comp_{c_0}(Y_1^r\cap\partial\MC)$ is the connected component of $Y_1^r\cap\partial\MC$ containing $c_0$. Since $\partial\MC$ is compact, there are finitely many $\Comp_{c_i}(Y_1^r(c_i)\cap\partial\MC)$'s such that their union covers $\partial\MC$. This implies that $\partial\MC$ is a quasi-circle.
\end{proof}

\begin{rmk}
If $\partial \MC\cap\Gamma\neq\emptyset$. By Corollary \ref{cor:singleton}, $\partial \MC\cap\Gamma$ is a singleton. The boundary $\partial\MC\setminus\Gamma$ is locally a quasi-arc by Theorem \ref{thm:quasi-disjoint}. To obtain that $\partial\MC$ is a quasi-circle in this case, one first needs to prove that the point $\partial \MC\cap\Gamma$ is not a cusp of $\partial\MC$.
However we cannot use the holomorphic motion in any neighborhood of such a point.
\end{rmk}

\section{Number of capture components}\label{sec:num-cap}

We calculate the number of capture components in $\Gamma_\Ext$ of each level. By Lemma \ref{para-c} it is sufficient to count the number of the centers of the capture components.

\begin{lem}\label{lem:sim-root}
The roots of $P_c^{\circ \ell}(c)=0$ are simple, where $\ell\geqslant 1$.
\end{lem}

\begin{proof}
Let $c_0$ be a root of $P_c^{\circ \ell}(c)=0$. Suppose $c_0$ is not simple. Then in a neighborhood of $c_0$ one can write
\begin{equation}
P_c^{\circ \ell}(c)=b\,(c-c_0)^k+\MO((c-c_0)^{k+1})
\end{equation}
with $b\neq 0$ and $k\geqslant 2$. Since $c_0$ belongs to a capture component $\MC$ in $\Gamma_\Ext$ with level $1\leqslant m\leqslant \ell$, it implies that $(c,z)\mapsto\phi_c(z)$ is well defined and holomorphic in a neighborhood of $(c_0,0)$, where $\phi_c$ is the linearizing map satisfying \eqref{equ:linear-map}. Note that in a small neighborhood of $c_0$ we have $\phi_c(z)=a_1(c)z+\MO(z^2)$ with $a_1(c)\neq 0$. Therefore, in a neighborhood of $c_0$,
\begin{equation}
\phi_c(P_c^{\circ \ell}(c))=b\,a_1(c)(c-c_0)^k+\MO((c-c_0)^{k+1}).
\end{equation}
On the other hand, in a neighborhood of $c_0$, we have
\begin{equation}
\phi_c(P_c^{\circ \ell}(c))=e^{2\pi\ii\theta(\ell-m)}\phi_c(P_c^{\circ m}(c))=e^{2\pi\ii\theta(\ell-m)}\Phi(c),
\end{equation}
where $\Phi:\MC\to\D$ is the conformal isomorphism in Lemma \ref{para-c}. It is clear that we obtain a contradiction by comparing the above two equations.
\end{proof}

\begin{thm}\label{thm:num-cap}
There are exactly $3^{\ell-1}$ capture components in $\Gamma_\Ext$ for any level $\ell\geqslant 1$.
\end{thm}

\begin{proof}
For simplicity, we denote by $\lambda=e^{2\pi\ii\theta}$. Note that $Q_0(c):=P_c^{\circ 0}(c)=c$ and
\begin{equation}
Q_1(c):=P_c(c)=\lambda\,c(3-c)/6=c\, G_1(c),
\end{equation}
where $G_1(c)=\lambda(3-c)/6$ is a linear map. Hence $\Gamma_\Ext$ contains exactly one capture component of level $1$ whose center is $c=3$.
Inductively, for $\ell\geqslant 2$, suppose $Q_{\ell-1}(c):=P_c^{\circ (\ell-1)}(c)=c\,G_{\ell-1}(c)$, where $G_{\ell-1}(c)$ is a polynomial.
Then a direct calculation shows that
\begin{equation}
Q_\ell(c):=P_c^{\circ \ell}(c)=P_c(Q_{\ell-1}(c))=c\,G_{\ell}(c),
\end{equation}
where
\begin{equation}
G_\ell(c)=\lambda\,G_{\ell-1}(c)\left(1-\frac{c(1+c)}{2}\,G_{\ell-1}(c)+\frac{c}{3}\Big(G_{\ell-1}(c)\Big)^2\right).
\end{equation}
Therefore we have
\begin{equation}
\deg G_\ell=3\deg G_{\ell-1}+1 \text{\quad and hence\quad} \deg G_\ell=(3^\ell-1)/2.
\end{equation}

According to Lemma \ref{lem:sim-root}, all the roots of $G_{\ell}(c)=0$ are simple. Note that the roots of $G_{\ell-1}(c)=0$ are the centers of the capture components with level at most $\ell-1$. Therefore, the number $N_\ell$ of the capture components with level $\ell\geqslant 1$ is
\begin{equation}
N_\ell=\deg G_\ell-\deg G_{\ell-1}=2\deg G_{\ell-1}+1=3^{\ell-1}.
\end{equation}
The proof is complete.
\end{proof}


\bibliographystyle{amsalpha}
\bibliography{E:/Latex-model/Ref1}

\providecommand{\bysame}{\leavevmode\hbox to3em{\hrulefill}\thinspace}
\providecommand{\MR}{\relax\ifhmode\unskip\space\fi MR }
\providecommand{\MRhref}[2]{%
  \href{http://www.ams.org/mathscinet-getitem?mr=#1}{#2}
}
\providecommand{\href}[2]{#2}
\begin{thebibliography}{QRWY15}

\bibitem[BCOT21]{BCOT21}
A.~Blokh, A.~Ch{\'{e}}ritat, L.~Oversteegen, and V.~Timorin, \emph{Location of
  {S}iegel capture polynomials in parameter spaces}, Nonlinearity \textbf{34}
  (2021), no.~4, 2430--2453.

\bibitem[BH88]{BH88}
B.~Branner and J.~H. Hubbard, \emph{The iteration of cubic polynomials. {I}.
  {T}he global topology of parameter space}, Acta Math. \textbf{160} (1988),
  no.~3-4, 143--206.

\bibitem[BH01]{BH01}
X.~Buff and C.~Henriksen, \emph{Julia sets in parameter spaces}, Comm. Math.
  Phys. \textbf{220} (2001), no.~2, 333--375.

\bibitem[BKM10]{BKM10}
A.~Bonifant, J.~Kiwi, and J.~Milnor, \emph{Cubic polynomial maps with periodic
  critical orbit. {II}. {E}scape regions}, Conform. Geom. Dyn. \textbf{14}
  (2010), 68--112.

\bibitem[BOPT14]{BOPT14}
A.~Blokh, L.~Oversteegen, R.~Ptacek, and V.~Timorin, \emph{The main cubioid},
  Nonlinearity \textbf{27} (2014), no.~8, 1879--1897.

\bibitem[BOST22]{BOST22}
A.~Blokh, L.~Oversteegen, A.~Shepelevtseva, and V.~Timorin, \emph{Modeling core
  parts of {Z}akeri slices {I}}, Mosc. Math. J. \textbf{22} (2022), no.~2,
  265--294.

\bibitem[BOT22]{BOT22}
A.~Blokh, L.~Oversteegen, and V.~Timorin, \emph{Slices of the parameter space
  of cubic polynomials}, Trans. Amer. Math. Soc. \textbf{375} (2022), no.~8,
  5313--5359.

\bibitem[CG93]{CG93}
L.~Carleson and T.~W. Gamelin, \emph{Complex dynamics}, Universitext: Tracts in
  Mathematics, Springer-Verlag, New York, 1993.

\bibitem[Ch{\'{e}}20]{Che20B}
A.~Ch{\'{e}}ritat, \emph{On the size of {S}iegel disks with fixed multiplier
  for cubic polynomials}, arXiv: 2003.13337, 2020.

\bibitem[Che25]{Che25}
D.~Cheraghi, \emph{Topology of irrationally indifferent attractors}, Ann. Sci.
  \'{E}c. Norm. Sup\'{e}r. (4) \textbf{58} (2025), no.~6, 1319--1381.

\bibitem[CWY22]{CWY22}
J.~Cao, X.~Wang, and Y.~Yin, \emph{Boundaries of capture hyperbolic
  components}, arXiv: 2206.07462, 2022.

\bibitem[DH85]{DH8485}
A.~Douady and J.~H. Hubbard, \emph{\'{E}tude dynamique des polyn\^{o}mes
  complexes. {P}artie {I}, {II}}, Publications Math\'{e}matiques d'Orsay, vol.
  84-85, Universit\'{e} de Paris-Sud, Orsay, 1984-1985.

\bibitem[DL23]{DL23a}
D.~Dudko and M.~Lyubich, \emph{Local connectivity of the {M}andelbrot set at
  some satellite parameters of bounded type}, Geom. Funct. Anal. \textbf{33}
  (2023), no.~4, 912--1047.

\bibitem[DL26]{DL26}
\bysame, \emph{M{LC} at {F}eigenbaum points}, Publ. Math. Inst. Hautes
  \'{E}tudes Sci. \textbf{143} (2026), 97--141.

\bibitem[Dou87]{Dou87}
A.~Douady, \emph{Disques de {S}iegel et anneaux de {H}erman}, Bourbaki seminar,
  {V}ol. 1986/87, Ast\'{e}risque, no. 152-153, Soc. Math. France, Paris, 1987,
  pp.~151--172.

\bibitem[EY99]{EY99}
A.~Epstein and M.~Yampolsky, \emph{Geography of the cubic connectedness locus:
  intertwining surgery}, Ann. Sci. \'{E}cole Norm. Sup. (4) \textbf{32} (1999),
  no.~2, 151--185.

\bibitem[Fau92]{Fau92}
D.~Faught, \emph{Local connectivity in a family of cubic polynomials}, Thesis
  (Ph.D.)--Cornell University, 1992.

\bibitem[FYZ22]{FYZ22}
Y.~Fu, F.~Yang, and G.~Zhang, \emph{Quadratic rational maps with a 2-cycle of
  {S}iegel disks}, J. Geom. Anal. \textbf{32} (2022), no.~10, Paper No. 244, 26
  pp.

\bibitem[GJW10]{GJW10}
F.~P. Gardiner, Y.~Jiang, and Z.~Wang, \emph{Holomorphic motions and related
  topics}, Geometry of {R}iemann surfaces, London Math. Soc. Lecture Note Ser.,
  vol. 368, Cambridge Univ. Press, Cambridge, 2010, pp.~156--193.

\bibitem[G{\'{S}}03]{GS03}
J.~Graczyk and G.~{\'{S}}wi{\c{a}}tek, \emph{Siegel disks with critical points
  in their boundaries}, Duke Math. J. \textbf{119} (2003), no.~1, 189--196.

\bibitem[Her87]{Her87}
M.~R. Herman, \emph{Conjugaison quasi-symm\'{e}trique des hom\'{e}omorphismes
  analytiques du cercle \`{a} des rotations},
  \href{https://www.math.univ-toulouse.fr/~cheritat/Herman/Herman_man_1_eng.pdf}{preliminary
  manuscript}, 1987.

\bibitem[Lav89]{Lav89}
P.~Lavaurs, \emph{Syst\`{e}mes dynamiques holomorphes: explosion de points
  periodiques paraboliques}, Thesis (Ph.D.)--Universit\'{e} de Paris-Sud,
  Orsay, 1989.

\bibitem[Lom14]{Lom14}
L.~Lomonaco, \emph{Parameter space for families of parabolic-like mappings},
  Adv. Math. \textbf{261} (2014), 200--219.

\bibitem[Lyu83]{Lyu83b}
M.~Lyubich, \emph{Some typical properties of the dynamics of rational
  mappings}, Uspekhi Mat. Nauk \textbf{38} (1983), no.~5(233), 197--198.

\bibitem[McM00]{McM00b}
C.~T. McMullen, \emph{The {M}andelbrot set is universal}, The {M}andelbrot set,
  theme and variations, London Math. Soc. Lecture Note Ser., vol. 274,
  Cambridge Univ. Press, Cambridge, 2000, pp.~1--17.

\bibitem[Mil92]{Mil92}
J.~Milnor, \emph{Remarks on iterated cubic maps}, Experiment. Math. \textbf{1}
  (1992), no.~1, 5--24.

\bibitem[Mil09]{Mil09}
\bysame, \emph{Cubic polynomial maps with periodic critical orbit. {I}},
  Complex dynamics, A K Peters, Wellesley, MA, 2009, pp.~333--411.

\bibitem[Mil12]{Mil12}
\bysame, \emph{Hyperbolic components, \emph{with an appendix by A. Poirier}},
  Conformal dynamics and hyperbolic geometry, Contemp. Math., vol. 573, Amer.
  Math. Soc., Providence, RI, 2012, pp.~183--232.

\bibitem[MSS83]{MSS83}
R.~Ma{\~{n}}{\'{e}}, P.~Sad, and D.~Sullivan, \emph{On the dynamics of rational
  maps}, Ann. Sci. \'{E}cole Norm. Sup. (4) \textbf{16} (1983), no.~2,
  193--217.

\bibitem[Nak05]{Nak05}
S.~Nakane, \emph{Capture components for cubic polynomials with parabolic fixed
  points}, Academic Reports Fac. Eng. Tokyo Polytech. Univ. \textbf{28} (2005),
  no.~1, 33--41.

\bibitem[PZ04]{PZ04}
C.~L. Petersen and S.~Zakeri, \emph{On the {J}ulia set of a typical quadratic
  polynomial with a {S}iegel disk}, Ann. of Math. (2) \textbf{159} (2004),
  no.~1, 1--52.

\bibitem[QRW23]{QRW23}
W.~Qiu, P.~Roesch, and Y.~Wang, \emph{Escape components of {M}c{M}ullen maps},
  Ergodic Theory Dynam. Systems \textbf{43} (2023), no.~11, 3745--3775.

\bibitem[QRWY15]{QRWY15}
W.~Qiu, P.~Roesch, X.~Wang, and Y.~Yin, \emph{Hyperbolic components of
  {M}c{M}ullen maps}, Ann. Sci. \'{E}c. Norm. Sup\'{e}r. (4) \textbf{48}
  (2015), no.~3, 703--737.

\bibitem[Roe07]{Roe07}
P.~Roesch, \emph{Hyperbolic components of polynomials with a fixed critical
  point of maximal order}, Ann. Sci. \'{E}cole Norm. Sup. (4) \textbf{40}
  (2007), no.~6, 901--949.

\bibitem[RWY17]{RWY17}
P.~Roesch, X.~Wang, and Y.~Yin, \emph{Moduli space of cubic {N}ewton maps},
  Adv. Math. \textbf{322} (2017), 1--59.

\bibitem[Shi98]{Shi98}
M.~Shishikura, \emph{The {H}ausdorff dimension of the boundary of the
  {M}andelbrot set and {J}ulia sets}, Ann. of Math. (2) \textbf{147} (1998),
  no.~2, 225--267.

\bibitem[Shi01]{Shi01}
\bysame, \emph{Herman's theorem on quasisymmetric linearization of analytic
  circle homemorphisms}, manuscript, 2001.

\bibitem[Sie42]{Sie42}
C.~L. Siegel, \emph{Iteration of analytic functions}, Ann. of Math. (2)
  \textbf{43} (1942), 607--612.

\bibitem[Slo91]{Slo91}
Z.~Slodkowski, \emph{Holomorphic motions and polynomial hulls}, Proc. Amer.
  Math. Soc. \textbf{111} (1991), no.~2, 347--355.

\bibitem[SY25]{SY25}
M.~Shishikura and F.~Yang, \emph{The high type quadratic {S}iegel disks are
  {J}ordan domains}, J. Eur. Math. Soc. (JEMS) \textbf{27} (2025), no.~11,
  4501--4562.

\bibitem[Wan21]{Wan21}
X.~Wang, \emph{Hyperbolic components and cubic polynomials}, Adv. Math.
  \textbf{379} (2021), No. 107554, 42 pp.

\bibitem[WY17]{WY17}
X.~Wang and Y.~Yin, \emph{Global topology of hyperbolic components: {C}antor
  circle case}, Proc. Lond. Math. Soc. (3) \textbf{115} (2017), no.~4,
  897--923.

\bibitem[YZ24]{YZ24}
J.~Yang and R.~Zhang, \emph{Rigidity of bounded type cubic {S}iegel
  polynomials}, arXiv: 2311.00431v2, 2024.

\bibitem[Zak99]{Zak99}
S.~Zakeri, \emph{Dynamics of cubic {S}iegel polynomials}, Comm. Math. Phys.
  \textbf{206} (1999), no.~1, 185--233.

\bibitem[Zak16]{Zak16}
\bysame, \emph{Conformal fitness and uniformization of holomorphically moving
  disks}, Trans. Amer. Math. Soc. \textbf{368} (2016), no.~2, 1023--1049.

\bibitem[Zak18]{Zak18}
\bysame, \emph{Rotation sets and complex dynamics}, Lecture Notes in
  Mathematics, vol. 2214, Springer, Cham, 2018.

\bibitem[Zha08]{Zha08b}
G.~Zhang, \emph{Dynamics of {S}iegel rational maps with prescribed
  combinatorics}, arXiv: 0811.3043, 2008.

\bibitem[Zha11]{Zha11}
\bysame, \emph{All bounded type {S}iegel disks of rational maps are
  quasi-disks}, Invent. Math. \textbf{185} (2011), no.~2, 421--466.

\bibitem[Zha24]{Zha24R}
R.~Zhang, \emph{On dynamical parameter space of cubic polynomials with a
  parabolic fixed point}, J. Lond. Math. Soc. (2) \textbf{110} (2024), no.~6,
  Paper No. e70038, 83 pp.

\end{thebibliography}

\end{document}